\documentclass[12pt,reqno,centertags,a4paper]{iopart}
\expandafter\let\csname equation*\endcsname\relax
\expandafter\let\csname endequation*\endcsname\relax

\usepackage[T1]{fontenc}
\usepackage{
amsfonts,
amssymb,
amsthm
}

\usepackage{pslatex}
\usepackage{enumerate}
\usepackage{nicefrac,xfrac}
\usepackage{latexsym}
\usepackage[classicReIm,uprightRoman,upright,frenchstyle,
uprightgreeks,notextcomp]{kpfonts}

\allowdisplaybreaks
\numberwithin{equation}{section}

\newcommand{\R}{\mathbb{R}}
\newcommand{\C}{\mathbb{C}}
\newcommand{\Z}{\mathbb{Z}}
\newcommand{\abs}[1]{\vert#1\vert}
\newcommand{\ol}[1]{\overline{#1}}
\newcommand{\mrm}[1]{\mathrm{#1}}
\newcommand{\co}{\!\colon\thinspace}
\newcommand{\Wh}[5]{ \left(\begin{array}{rr}#1, & #3 \\ #2, & #4 \end{array};#5\right) }

\begin{document}

\title{On the definite integral of two confluent hypergeometric functions related to 
the Kamp\'{e} de F\'{e}riet double series }

\author{ R. Jur\v{s}\.{e}nas }
 
\address{ Institute of Theoretical Physics and Astronomy of Vilnius University, \\
A.~Go\v{s}tauto 12, 01108 Vilnius, Lithuania }

\begin{abstract}
The Kamp\'{e} de F\'{e}riet double series $F_{1:1;1}^{1:1;1}$ is studied through the solution to the associated first-order 
nonhomogeneous differential equation. It is shown that the integral of 
$t^{\beta+l}M(\cdot;\beta;\lambda t)M(\cdot;\beta;-\lambda t)$ over $t\in[0,T]$, $T\geq0$, $l=0,1,\ldots$,
$\Re\beta+l>-1$, is a linear combination of functions $F_{1:1;1}^{1:1;1}$. The integral is a generalization of 
a class of so-called Coulomb integrals involving regular Coulomb wave functions.
\end{abstract}


{\small
\tableofcontents
}

\maketitle

\section{Introduction}\label{sec:s1}

In this paper we study the Kamp\'{e} de F\'{e}riet function $F_{1:1;1}^{1:1;1}$ as a solution to the first-order
nonhomogeneous differential equation on $\C$,

\begin{equation}
\left(-\frac{d}{dz}+\frac{\beta}{z}\right)f(z)=\frac{ \Gamma(\alpha+1)\Gamma(\alpha+\gamma) }{ z\Gamma(\beta)\Gamma(\delta) }
M(\alpha+1;\delta;-1/z)M(\alpha+\gamma;\beta;1/z)
\label{eq:diffeq}
\end{equation}

\noindent{}where $\Gamma$ is the gamma function, $M$ is the Kummer's function \cite[\S 6, \S 13]{Abr72}.
Throughout, the parameters $\alpha$, $\beta$, $\gamma$, $\delta\in\C$ and the argument $z\in\C$ are assumed to be such that the 
hypergeometric functions as well as the gamma functions are well-defined unless additional conditions are imposed.

\subsection{Main results : Solution}

The solution to (\ref{eq:diffeq}) reads

\begin{equation}
f(z)\equiv\Wh{\alpha}{\beta}{\gamma}{\delta}{z}=
\frac{\Gamma(\alpha+1)\Gamma(\alpha+\gamma)}{\Gamma(\beta+1)\Gamma(\delta)}\:
F_{1:1;1}^{1:1;1}\left(\begin{array}{rrr}\beta: & \alpha+1; & \alpha+\gamma; \\
\beta+1: & \delta; & \beta; \end{array}-\frac{1}{z},\frac{1}{z}\right)
\label{eq:hyper3}
\end{equation}

\noindent{}plus the solution $Cz^{\beta}$ to the homogeneous equation. As long as the constant of integration $C\in\C$ is arbitrary,
we choose $C=0$.

\subsection{Main results : Definite integral}

The integral of interest is the following,

\begin{align}
&\int_{0}^{T}t^{\beta+l}M(\alpha;\beta;\lambda t)M(\alpha+\gamma;\beta;-\lambda t)dt=
\frac{ \Gamma(\beta)\Gamma(\beta+l+1)T^{\beta+l+1} }{ \Gamma(\alpha+\gamma) }\sum_{n=0}^{l+1}
\frac{ (-\lambda T)^{n}\Gamma(\alpha+\gamma+n) }{ \Gamma(\beta+n)\Gamma(\beta+l+n+2) } \nonumber \\
&\times\binom{l+1}{n}
F_{1:1;1}^{1:1;1}\left(\begin{array}{rrr}\beta+l+n+1: & \alpha; & \alpha+\gamma+n; \\
\beta+l+n+2: & \beta; & \beta+l+n+1; \end{array}\lambda T,-\lambda T\right) \nonumber \\
&(\Re\beta+l>-1,\quad l=0,1,\ldots,\quad \lambda\in\C,\quad T\geq0).
\label{eq:integralMM}
\end{align}

\noindent{}The Barnes integral representation is shown in (\ref{eq:tl-2}).

\subsection{Novelty and connection with known results}

The class of integrals (\ref{eq:integralMM}), with $\beta\in\Z_{+}$ and $\beta+l\leq2$, has been considered in
\cite{Arnoldus92} when dealing with the Coulomb wave functions $F_{L}(\eta,\rho)$ ($L=0,1,\ldots$ is the
angular momentum; $\eta\in\R\backslash\{0\}$ is the coupling constant; $\rho\geq0$ is the distance) \cite[\S 14]{Abr72}.
It follows that the right-hand side of (\ref{eq:integralMM}) extends the results for $-1<\Re\beta+l\leq2$ as well as brings forth
new results for $\Re\beta+l>2$. In particular, (\ref{eq:integralMM}) yields

\begin{align}
\int_{0}^{T}\rho^{l}\;\ol{ F_{L}(\eta,\rho) }F_{L}(\eta,\rho)d\rho=&
\frac{ 4^{L}e^{-\pi\eta}\Gamma(L+1-i\eta)\Gamma(2L+l+3)T^{2L+l+3} }{ \Gamma(2L+2) } \nonumber \\
&\times \sum_{n=0}^{l+1}\frac{ (-2iT)^{n}\Gamma(L+n+1+i\eta) }{ \Gamma(2L+n+2)\Gamma(2L+l+n+4) }
\binom{l+1}{n} \nonumber \\
&\times 
F_{1:1;1}^{1:1;1}\left(\begin{array}{rrr}2L+l+n+3: & L+1-i\eta; & L+n+1+i\eta; \\
2L+l+n+4: & 2L+2; & 2L+l+n+3; \end{array}2iT,-2iT\right) \nonumber \\
&(L,l\in\Z_{+}\cup\{0\},\quad \eta\in\R\backslash\{0\},\quad T\geq0)
\label{eq:CoulombIntegral}
\end{align}

\noindent{}where $^{-}\co\C\to\C$ is the complex conjugation; $i\equiv\sqrt{-1}$.

Unlike the integrals considered by other authors (see eg \cite{Arnoldus92,Gradshteyn07,Hamza74,Saad03,Shanker68} and the citations
therein), the integrals (\ref{eq:integralMM})--(\ref{eq:CoulombIntegral}) diverge as $T\to\infty$, and most likely, this
is the main reason why they were not so widely used in various applications. However, recent results on
the particles inside a finite cubic box \cite{Kreuzer10,Kreuzer12} add up the present class of integrals to those
with possible request.

In the present paper, we also study the asymptotic expansion (\S\S\ref{sec:conv}--\ref{sec:asym}) of solution 
(\ref{eq:hyper3}) as well as derive some recurrence relations (\S\ref{sec:recur}). In virtue of (\ref{eq:integralMM}),
we give an alternative proof (\S\ref{sec:exp}) of the integral

\begin{equation}
\int_{0}^{\infty}t^{\beta}e^{-ht}M(\alpha;\beta;it)M(\alpha+\gamma;\beta;-it)dt=
\Gamma(\beta+1)h^{-\beta-1}\:
F_{2}(\beta+1;\alpha,\alpha+\gamma;\beta,\beta;i/h,-i/h)
\label{eq:F2}
\end{equation}

\noindent{}($\Re\beta>-1;\vert h\vert>2$) studied in \cite{Saad03};
$F_{2}(\cdot)$ denotes the Appell's hypergeometric function (see eg \cite{Vidunas09}).

\subsection{Notations}

The Kamp\'e de F\'eriet function $F_{1:1;1}^{1:1;1}$ of two variables is defined in agreement with \cite{Cvi10,Gro91},

\begin{equation}
F_{1:1;1}^{1:1;1}\left(\begin{array}{rrr}a_{1}: & b_{1}; & b_{2}; \\
c_{1}: & d_{1}; & d_{2}; \end{array}z_{1},z_{2}\right)=\sum_{m_{1},m_{2}=0}^{\infty}
\frac{(a_{1})_{m_{1}+m_{2}} (b_{1})_{m_{1}} (b_{2})_{m_{2}} z_{1}^{m_{1}}z_{2}^{m_{2}} }{
(c_{1})_{m_{1}+m_{2}} (d_{1})_{m_{1}} (d_{2})_{m_{2}} m_{1}!m_{2}! }
\label{eq:Feriet}
\end{equation}

\noindent{}($(a)_{m}$ the Pochhammer symbol).

The Kummer's function $M(\cdot)\equiv\,_{1}F_{1}(\cdot)$, the generalized hypergeometric function $_{p}F_{q}(\cdot)$,
and the MacRobert's $E$-function are related to each other through the relations \cite{Mac54,Mac58}

\begin{subequations}\label{eq:MacRob}
\begin{align}
E\left(\begin{array}{rrrr}a_{1}, & a_{2}, & \ldots, & a_{p} \\ b_{1}, & b_{2}, & \ldots, & b_{q} \end{array};z\right)=&
\frac{\prod_{j=1}^{p}\Gamma(a_{j})}{\prod_{j=1}^{q}\Gamma(b_{j})}\:
_{p}F_{q}\left(\begin{array}{rrrr}
a_{1}, & a_{2}, & \ldots, & a_{p} \\ b_{1}, & b_{2}, & \ldots, & b_{q} \end{array};-\frac{1}{z}\right)\;\;\;
(p\leq q)
\label{eq:MacRob-1} \\
=&\frac{1}{2\pi i}\int_{B}\frac{z^{\zeta}\Gamma(\zeta)\prod_{j=1}^{p}\Gamma(a_{j}-\zeta)}{\prod_{j=1}^{q}\Gamma(b_{j}-\zeta)}d\zeta
\;\;\;(\vert\mrm{arg}\;z\vert<\pi) \label{eq:MacRob-2}
\end{align}
\end{subequations}

\noindent{}with suitable $\{a_{j}\in\C\co j=1,\ldots,p\}$, $\{b_{j}\in\C\co j=1,\ldots,q\}$, $z\in\C$.
For $p=q=1$ one writes $E(a_{1};b_{1};z)$.
The Barnes-type contour $B$ is taken up the $\Im\zeta$-axis so that the poles of $\{\Gamma(a_{j}-\zeta)\co j$ $=1,\ldots,p\}$ 
and $\{\Gamma(b_{j}-\zeta)\co j=1,\ldots,q\}$ lie to the right of the contour and the poles of $\Gamma(\zeta)$ to the left of the 
contour.

The notation of $E$-function is convenient for the proof of (\ref{eq:integralMM}) as well as for simplification of the 
recurrence relations \S\ref{sec:recur}.

\section{Proof of Eq.~(\ref{eq:hyper3})}

It suffices to prove that the solution to (\ref{eq:diffeq}) can be represented as a series

\begin{equation}
f(z)\equiv\Wh{\alpha}{\beta}{\gamma}{\delta}{z}=\frac{\Gamma(\alpha+1)}{\Gamma(\delta)}
\sum_{\nu=0}^{\infty}\frac{z^{-\nu}\Gamma(\alpha+\gamma+\nu)}{\nu!\Gamma(\beta+1+\nu)}\:
_{2}F_{2}\left(\begin{array}{rr}\alpha+1, & \beta+\nu \\ \delta, & \beta+1+\nu \end{array};-\frac{1}{z}\right).
\label{eq:hyper}
\end{equation}

\noindent{}Then the proof of (\ref{eq:hyper3}) is straightforward: Rewrite $_{2}F_{2}(\cdot)$ in (\ref{eq:hyper})
as an infinite series; then (\ref{eq:hyper3}) holds by virtue of (\ref{eq:Feriet}).

In order to prove (\ref{eq:hyper}), make use of the relation

\begin{equation}
z\frac{d}{dz}\:_{2}F_{2}\left(\begin{array}{rr}\alpha+1, & \beta+\nu \\ \delta, & \beta+1+\nu \end{array};z\right)=
(\beta+\nu)\left(M(\alpha+1;\delta;z)-\:
_{2}F_{2}\left(\begin{array}{rr}\alpha+1, & \beta+\nu \\ \delta, & \beta+1+\nu \end{array};z\right)
\right).
\label{eq:diffE}
\end{equation}

\noindent{}Since $z\to-1/z$ implies $z(d/dz)\to-z(d/dz)$, assign $z(d/dz)$ to (\ref{eq:hyper}), then substitute 
(\ref{eq:diffE}) in obtained expression and get that

\begin{align*}
\frac{\Gamma(\delta)}{\Gamma(\alpha+1)}z\frac{d}{dz}\Wh{\alpha}{\beta}{\gamma}{\delta}{z}=&
-\sum_{\nu=1}^{\infty}\frac{z^{-\nu}\Gamma(\alpha+\gamma+\nu)}{(\nu-1)!\Gamma(\beta+1+\nu)}\:
_{2}F_{2}\left(\begin{array}{rr}\alpha+1, & \beta+\nu \\ \delta, & \beta+1+\nu \end{array};-\frac{1}{z}\right) \\
&+\sum_{\nu=0}^{\infty}\frac{z^{-\nu}\Gamma(\alpha+\gamma+\nu)}{\nu!\Gamma(\beta+1+\nu)}\:z\frac{d}{dz}\:
_{2}F_{2}\left(\begin{array}{rr}\alpha+1, & \beta+\nu \\ \delta, & \beta+1+\nu \end{array};-\frac{1}{z}\right) \\
=&-\frac{\Gamma(\delta)}{z\Gamma(\alpha+1)}\Wh{\alpha}{\beta+1}{\gamma+1}{\delta}{z}+
\sum_{\nu=0}^{\infty}\frac{z^{-\nu}\Gamma(\alpha+\gamma+\nu)}{\nu!\Gamma(\beta+\nu)} \\
&\times\left(-M(\alpha+1;\delta;-1/z)+\:
_{2}F_{2}\left(\begin{array}{rr}\alpha+1, & \beta+\nu \\ \delta, & \beta+1+\nu \end{array};-\frac{1}{z}\right) \right) \\
=&-\frac{\Gamma(\delta)}{z\Gamma(\alpha+1)}\Wh{\alpha}{\beta+1}{\gamma+1}{\delta}{z}-
\frac{\Gamma(\alpha+\gamma)}{\Gamma(\beta)}M(\alpha+1;\delta;-1/z) \\
&\times M(\alpha+\gamma;\beta;1/z)
+\frac{\beta\Gamma(\delta)}{\Gamma(\alpha+1)}\Wh{\alpha}{\beta}{\gamma}{\delta}{z} \\
&+\sum_{\nu=1}^{\infty}\frac{z^{-\nu}\Gamma(\alpha+\gamma+\nu)}{(\nu-1)!\Gamma(\beta+1+\nu)}
\:_{2}F_{2}\left(\begin{array}{rr}\alpha+1, & \beta+\nu \\ \delta, & \beta+1+\nu \end{array};-\frac{1}{z}\right) \\
=&-\frac{\Gamma(\alpha+\gamma)}{\Gamma(\beta)}M(\alpha+1;\delta;-1/z)M(\alpha+\gamma;\beta;1/z) \\
&+\frac{\beta\Gamma(\delta)}{\Gamma(\alpha+1)}\Wh{\alpha}{\beta}{\gamma}{\delta}{z}\quad 
\text{as required}.
\end{align*}

\section{Proof of Eq.~(\ref{eq:integralMM})}

To simplify algebraic manipulations, we found $E$-functions more convenient compared to the
generalized hypergeometric functions. For this, let us consider the integral

\begin{equation}
\int_{0}^{T}t^{\beta+l}E(\alpha;\beta;\lambda/t)E(\alpha+\gamma;\beta;-\lambda/t)dt
\label{eq:integralEE}
\end{equation}

\noindent{}with $l\in\Z_{+}\cup\{0\}$ and suitable $\lambda\in\C$. Then (\ref{eq:integralMM}) is easy to
recover from (\ref{eq:integralEE}) by making the substitution $\lambda\to-1/\lambda$ in virtue of (\ref{eq:MacRob-1}).

By using the series representation of $E$-function, the integral (\ref{eq:integralEE}) can be rewritten as follows

\begin{equation}
\int_{0}^{T}t^{\beta+l}E(\alpha;\beta;\lambda/t)E(\alpha+\gamma;\beta;-\lambda/t)dt=
\sum_{\mu,\nu=0}^{\infty}
\frac{(-1/\lambda)^{\mu}(1/\lambda)^{\nu}\Gamma(\alpha+\mu)\Gamma(\alpha+\gamma+\nu)}{\mu!\nu!\Gamma(\beta+\mu)\Gamma(\beta+\nu)}
\int_{0}^{T}t^{\beta+l+\mu+\nu}dt. 
\label{eq:tdt}
\end{equation}

\noindent{}The integral over $t\in[0,T]$ exists and is equal to
$T^{\beta+l+\mu+\nu+1}/(\beta+l+\mu+\nu+1)$ provided $\Re\beta+l>-1$.
Expand $(\beta+l+\mu+\nu+1)^{-1}$ into the binomial series,

\begin{equation}
\frac{1}{\beta+l+\mu+\nu+1}=\frac{1}{\beta+1}\sum_{r,s=0}^{\infty}\frac{(-1)^{r+s}\nu^{s}(r+s)!}{r!s!(\beta+1)^{r+s}}
\sum_{u=0}^{r}\binom{r}{u}\mu^{u}l^{r-u}.
\label{eq:denom}
\end{equation}

\noindent{}Substitute (\ref{eq:denom}) in (\ref{eq:tdt}) and get that

\begin{align}
&\int_{0}^{T}t^{\beta+l}E(\alpha;\beta;\lambda/t)E(\alpha+\gamma;\beta;-\lambda/t)dt=\frac{T^{\beta+l+1}}{\beta+1}
\sum_{r,s=0}^{\infty}\frac{(-1)^{r+s}(r+s)!}{r!s!(\beta+1)^{r+s}} \nonumber \\
&\times \left(\sum_{u=0}^{r}\binom{r}{u}l^{r-u}
\sum_{\mu=0}^{\infty}\frac{(-T/\lambda)^{\mu}\mu^{u}\Gamma(\alpha+\mu)}{\mu!\Gamma(\beta+\mu)}\right)
\sum_{\nu=0}^{\infty}\frac{(T/\lambda)^{\nu}\nu^{s}\Gamma(\alpha+\gamma+\nu)}{\nu!\Gamma(\beta+\nu)}. \label{eq:int2}
\end{align}

\noindent{}The sum over $\nu=0,1,\ldots$ is equal to

\begin{enumerate}[\upshape (a)]
\item $s=0$: $E(\alpha+\gamma;\beta;-\lambda/T)$
\item $s=1$: $(T/\lambda)E(\alpha+\gamma+1;\beta+1;-\lambda/T)$
\item $s=2,3,\ldots$: 
\begin{align*}
\sum_{\nu=0}^{\infty}\frac{(T/\lambda)^{\nu}\nu^{s}\Gamma(\alpha+\gamma+\nu)}{\nu!\Gamma(\beta+\nu)}=&
(T/\lambda)
\sum_{\nu=0}^{\infty}\frac{(T/\lambda)^{\nu}(1+\nu)^{s-1}\Gamma(\alpha+\gamma+\nu+1)}{\nu!\Gamma(\beta+\nu+1)} \\
=&(T/\lambda)
E\left(\begin{array}{rr}\alpha+\gamma+1, & [2]_{s-1} \\ \beta+1, & [1]_{s-1} \end{array};-\frac{\lambda}{T}\right)
\end{align*}
\end{enumerate}

\noindent{}The symbol $[a]_{n}$ ($a\in\C$; $n=1,2,\ldots$) denotes $n$ copies of $a$,

\begin{equation}
E\left(\begin{array}{rr}\alpha+\gamma+1, & [2]_{n} \\ \beta+1, & [1]_{n} \end{array};-\frac{\lambda}{T}\right)\equiv
E\left(\begin{array}{rrrrr}\alpha+\gamma+1, & 2, & 2, & \ldots & 2,
\\ \beta+1, & 1, & 1, & \ldots & 1, \end{array};-\frac{\lambda}{T}\right);
\label{eq:defnE}
\end{equation}

\noindent{}we set $[a]_{0}$ as being empty so that $\text{(c)}\Rightarrow\text{(b)}$ for $s=1$.

The sum over $\mu=0,1,\ldots$ in (\ref{eq:int2}) is found from (a)--(c) with the substitutions $\lambda\to-\lambda$,
$s\to u$, $\alpha+\gamma\to\alpha$.

The sum over $r=0,1,\ldots$ is extracted into the term with $r=0$ and those with $r=1,2,\ldots$ For $r=0$, the sum
in the parentheses in (\ref{eq:int2}) is given by $E(\alpha;\beta;\lambda/T)$. For $r=1,2,\ldots$, it is equal to

\begin{align}
&l^{r}E(\alpha;\beta;\lambda/T)-(T/\lambda)\sum_{u=1}^{r}\binom{r}{u}l^{r-u}
E\left(\begin{array}{rr}\alpha+1, & [2]_{u-1} \\ \beta+1, & [1]_{u-1} \end{array};
\frac{\lambda}{T}\right)=l^{r}E(\alpha;\beta;\lambda/T) \nonumber \\
&-(T/\lambda)\sum_{\nu=0}^{\infty}\frac{(-T/\lambda)^{\nu}\Gamma(\alpha+1+\nu)}{(\nu+1)!\Gamma(\beta+1+\nu)}\sum_{u=1}^{r}
\binom{r}{u}(\nu+1)^{u}l^{r-u}=
l^{r}E(\alpha;\beta;\lambda/T) \nonumber \\
&+(T/\lambda)l^{r}\sum_{\nu=0}^{\infty}\frac{(-T/\lambda)^{\nu}\Gamma(\alpha+1+\nu)}{(\nu+1)!\Gamma(\beta+1+\nu)}-(T/\lambda)
\sum_{\nu=0}^{\infty}\frac{(-T/\lambda)^{\nu}\Gamma(\alpha+1+\nu)(l+1+\nu)^{r}}{(\nu+1)!\Gamma(\beta+1+\nu)} \nonumber \\
&=l^{r}E(\alpha;\beta;\lambda/T)+l^{r}\left(\frac{\Gamma(\alpha)}{\Gamma(\beta)}-E(\alpha;\beta;\lambda/T)\right)+
\sum_{\nu=1}^{\infty}\frac{(-T/\lambda)^{\nu}\Gamma(\alpha+\nu)(l+\nu)^{r}}{\nu!\Gamma(\beta+\nu)} \nonumber \\
&=E\left(\begin{array}{rr}\alpha, & [l+1]_{r} \\ \beta, & [l]_{r} \end{array};\frac{\lambda}{T}\right).
\label{eq:sumr}
\end{align}

\noindent{}[Note that (\ref{eq:sumr}) yields $E(\alpha;\beta;\lambda/T)$ which corresponds to $r=0$.]
Substitute (a)--(c) and (\ref{eq:sumr}) in (\ref{eq:int2}). The result reads

\begin{align}
&\int_{0}^{T}t^{\beta+l}E(\alpha;\beta;\lambda/t)E(\alpha+\gamma;\beta;-\lambda/t)dt=\frac{T^{\beta+l+1}}{\beta+1}\Biggl[
E(\alpha;\beta;\lambda/T)E(\alpha+\gamma;\beta;-\lambda/T) \nonumber \\
&+(T/\lambda)E(\alpha;\beta;\lambda/T)\sum_{s=1}^{\infty}\frac{(-1)^{s}}{(\beta+1)^{s}}\:
E\left(\begin{array}{rr}\alpha+\gamma+1, & [2]_{s-1} \\ \beta+1, & [1]_{s-1} \end{array};-\frac{\lambda}{T}\right) \nonumber \\
&+E(\alpha+\gamma;\beta;-\lambda/T)\sum_{s=1}^{\infty}\frac{(-1)^{s}}{(\beta+1)^{s}}\:
E\left(\begin{array}{rr}\alpha, & [l+1]_{s} \\ \beta, & [l]_{s} \end{array};\frac{\lambda}{T}\right)+(T/\lambda)\sum_{r=1}^{\infty}
\frac{(-1)^{r}}{r!(\beta+1)^{r}} \nonumber \\
&\times E\left(\begin{array}{rr}\alpha, & [l+1]_{r} \\ \beta, & [l]_{r} \end{array};\frac{\lambda}{T}\right)
\sum_{s=1}^{\infty}\frac{(-1)^{s}(r+s)!}{s!(\beta+1)^{s}}\:
E\left(\begin{array}{rr}\alpha+\gamma+1, & [2]_{s-1} \\ \beta+1, & [1]_{s-1} \end{array};-\frac{\lambda}{T}\right) \Biggr].
\label{eq:sumr2}
\end{align}

Next, we prove that

\begin{subequations}\label{eq:Sum}
\begin{align}
&\sum_{s=1}^{\infty}\frac{(-1)^{s}}{(\beta+1)^{s}}\:E\left(\begin{array}{rr}\alpha+1, & [2]_{s-1} \\
\beta+1, & [1]_{s-1} \end{array};z\right)=-E\left(\begin{array}{rr}\alpha+1, & \beta+2 \\ \beta+1, & \beta+3 \end{array};z\right),
\label{eq:Sum1} \\
&\sum_{s=1}^{\infty}\frac{(-1)^{s}(r+s)!}{s!(\beta+1)^{s}}\:
E\left(\begin{array}{rr}\alpha+1, & [2]_{s-1} \\ \beta+1, & [1]_{s-1} \end{array};z\right)=zr!
\Biggl(E(\alpha;\beta;z) \nonumber \\
&-(\beta+1)^{r+1}E\left(\begin{array}{rr}\alpha, & [\beta+1]_{r+1} \\ \beta, & [\beta+2]_{r+1} \end{array};z\right)\Biggr)\quad
(r=0,1,\ldots), \label{eq:Sum2}.
\end{align}
\end{subequations}

Apply (\ref{eq:MacRob-2}) to the left-hand side of (\ref{eq:Sum1}) and write

\begin{align*}
&\sum_{s=1}^{\infty}\frac{(-1)^{s}}{(\beta+1)^{s}}\:E\left(\begin{array}{rr}\alpha+1, & [2]_{s-1} \\
\beta+1, & [1]_{s-1} \end{array};z\right)=\frac{1}{2\pi i}\int_{B}d\zeta\:
\frac{z^{\zeta}\Gamma(\zeta)\Gamma(\alpha+1-\zeta)}{\Gamma(\beta+1-\zeta)} \\
&\times\sum_{s=1}^{\infty}\frac{(-1)^{s}(1-\zeta)^{s-1}}{(\beta+1)^{s}}=-\frac{1}{2\pi i}\int_{B}d\zeta\:
\frac{z^{\zeta}\Gamma(\zeta)\Gamma(\alpha+1-\zeta)}{\Gamma(\beta+1-\zeta)}\cdot\frac{1}{\beta+2-\zeta}.
\end{align*}

\noindent{}By applying (\ref{eq:MacRob-2}) once again along with the substitution
$(\beta+2-\zeta)^{-1}=\Gamma(\beta+2-\zeta)/\Gamma(\beta+3-\zeta)$ one derives (\ref{eq:Sum1}).

Similarly, apply (\ref{eq:MacRob-2}) to the left-hand side of (\ref{eq:Sum2}) and get

\begin{align}
&\sum_{s=1}^{\infty}\frac{(-1)^{s}(r+s)!}{s!(\beta+1)^{s}}\:
E\left(\begin{array}{rr}\alpha+1, & [2]_{s-1} \\ \beta+1, & [1]_{s-1} \end{array};z\right)=\frac{1}{2\pi i}\int_{B}d\zeta\:
\frac{z^{\zeta}\Gamma(\zeta)\Gamma(\alpha+1-\zeta)}{\Gamma(\beta+1-\zeta)} \nonumber \\
&\times\sum_{s=1}^{\infty}\frac{(-1)^{s}(r+s)!(1-\zeta)^{s-1}}{s!(\beta+1)^{s}}=\frac{1}{2\pi i}\int_{B}d\zeta\:
\frac{z^{\zeta}\Gamma(\zeta)\Gamma(\alpha+1-\zeta)}{\Gamma(\beta+1-\zeta)}\cdot\frac{r!}{1-\zeta} \nonumber \\
&\times\left(-1+\frac{(\beta+1)^{r+1}}{(\beta+2-\zeta)^{r+1}} \right)=-\frac{r!}{2\pi i}\int_{B}d\zeta\:
\frac{z^{\zeta}\Gamma(\zeta)\Gamma(\alpha+1-\zeta)}{\Gamma(\beta+1-\zeta)(1-\zeta)}+\frac{r!(\beta+1)^{r+1}}{2\pi i} \nonumber \\
&\times \int_{B}d\zeta\:\frac{z^{\zeta}\Gamma(\zeta)\Gamma(\alpha+1-\zeta)}{\Gamma(\beta+1-\zeta)(1-\zeta)(\beta+2-\zeta)^{r+1}}.
\label{eq:Sum2a}
\end{align}

\noindent{}In order to calculate the obtained contour integrals, two possibilities are valid: either apply (\ref{eq:MacRob-2})
and then reduce the order of obtained $E$-functions or make use of the residues directly recalling that 
$\mrm{Res}_{\zeta=-\nu}\:\Gamma(\zeta)=(-1)^{\nu}/\nu!$ ($\nu=0,1,\ldots$). The combination of both strategies yields

\begin{align*}
&\frac{1}{2\pi i}\int_{B}d\zeta\:\frac{z^{\zeta}\Gamma(\zeta)\Gamma(\alpha+1-\zeta)}{\Gamma(\beta+1-\zeta)(1-\zeta)}=
E\left(\begin{array}{rr}\alpha+1, & 1 \\ \beta+1, & 2 \end{array};z\right)=z\left(\frac{\Gamma(\alpha)}{\Gamma(\beta)}-
E(\alpha;\beta;z)\right), \\
&\frac{1}{2\pi i}\int_{B}d\zeta\:
\frac{z^{\zeta}\Gamma(\zeta)\Gamma(\alpha+1-\zeta)}{\Gamma(\beta+1-\zeta)(1-\zeta)(\beta+2-\zeta)^{r+1}}=
E\left(\begin{array}{rrr}\alpha+1, & [\beta+2]_{r+1}, & 1 \\ \beta+1, & [\beta+3]_{r+1}, & 2 \end{array};z\right) \\
&=\sum_{\nu=0}^{\infty}\frac{(-1/z)^{\nu}\Gamma(\alpha+1+\nu)}{(\nu+1)!\Gamma(\beta+1+\nu)(\beta+2+\nu)^{r+1}}=
\frac{z\Gamma(\alpha)}{\Gamma(\beta)(\beta+1)^{r+1}} \\
&+\sum_{\nu=0}^{\infty}\frac{(-1/z)^{\nu-1}\Gamma(\alpha+\nu)}{ \nu!\Gamma(\beta+\nu)(\beta+1+\nu)^{r+1} }=
\frac{z\Gamma(\alpha)}{\Gamma(\beta)(\beta+1)^{r+1}}
-zE\left(\begin{array}{rr}\alpha, & [\beta+1]_{r+1} \\ \beta, & [\beta+2]_{r+1} \end{array};z\right).
\end{align*}

\noindent{}Substitute obtained expressions in (\ref{eq:Sum2a}) and get (\ref{eq:Sum2}).

By making proper replacements of the parameters substitute (\ref{eq:Sum1})--(\ref{eq:Sum2}) in (\ref{eq:sumr2}), make use of
the equation [in \cite[eq.~(4.4)]{Bho62} substitute $p=q=2$, $\alpha_{1}=\alpha+\gamma$, $\alpha_{2}=\beta+1$,
$\rho_{1}=\beta$, $\rho_{2}=\beta+2$, and after that apply \cite[\S5.2.1, eq.~(7)]{Erd53} with $a_{1}=\alpha+\gamma-1$, $\rho_{1}=
\beta-1$]

$$
(\beta+1)E\left(\begin{array}{rr}\alpha+\gamma, & \beta+1 \\ \beta, & \beta+2 \end{array};z\right)=
E(\alpha+\gamma;\beta;z)+(1/z)
E\left(\begin{array}{rr}\alpha+\gamma+1, & \beta+2 \\ \beta+1, & \beta+3 \end{array};z\right)
$$

\noindent{}and get that

\begin{align}
&\int_{0}^{T}t^{\beta+l}E(\alpha;\beta;\lambda/t)E(\alpha+\gamma;\beta;-\lambda/t)dt
=T^{\beta+l+1}\Biggl[E(\alpha;\beta;\lambda/T) 
E\left(\begin{array}{rr}\alpha+\gamma, & \beta+1 \\ \beta, & \beta+2 \end{array};-\frac{\lambda}{T}\right) \nonumber \\
&+\sum_{s=1}^{\infty}(-1)^{s}
E\left(\begin{array}{rr}\alpha, & [l+1]_{s} \\ \beta, & [l]_{s} \end{array};\frac{\lambda}{T}\right)
E\left(\begin{array}{rr}\alpha+\gamma, & [\beta+1]_{s+1} \\ \beta, & [\beta+2]_{s+1} \end{array};
-\frac{\lambda}{T}\right) \Biggr].
\label{eq:tl-1}
\end{align}

\noindent{}By (\ref{eq:MacRob-2}), the latter sum can be represented by

\begin{align*}
&\sum_{s=1}^{\infty}(-1)^{s}E\left(\begin{array}{rr}\alpha, & [l+1]_{s} \\ \beta, & [l]_{s} \end{array};\frac{\lambda}{T}\right)
E\left(\begin{array}{rr}\alpha+\gamma, & [\beta+1]_{s+1} \\ \beta, & [\beta+2]_{s+1} \end{array};-\frac{\lambda}{T}\right) \\
&=\frac{1}{2\pi i}\int_{B}d\zeta\:\frac{(\lambda/T)^{\zeta}\Gamma(\zeta)\Gamma(\alpha-\zeta)}{\Gamma(\beta-\zeta)}\cdot
\frac{1}{2\pi i}\int_{B}
d\zeta^{\prime}\:
\frac{(-\lambda/T)^{\zeta^{\prime}}\Gamma(\zeta^{\prime})\Gamma(\alpha+\gamma-\zeta^{\prime})}{\Gamma(\beta-\zeta^{\prime})} \\
&\times \sum_{s=1}^{\infty}\frac{(-1)^{s}(l-\zeta)^{s}}{(\beta+1-\zeta^{\prime})^{s+1}}=
\frac{1}{2\pi i}\int_{B}d\zeta\:\frac{(\lambda/T)^{\zeta}\Gamma(\zeta)\Gamma(\alpha-\zeta)}{\Gamma(\beta-\zeta)}\cdot
\frac{1}{2\pi i}\int_{B}
d\zeta^{\prime} \\
&\times \frac{(-\lambda/T)^{\zeta^{\prime}}\Gamma(\zeta^{\prime})\Gamma(\alpha+\gamma-\zeta^{\prime})}{\Gamma(\beta-\zeta^{\prime})}
\left(-\frac{1}{\beta+1-\zeta^{\prime}}+\frac{1}{\beta+l+1-\zeta-\zeta^{\prime}} \right) \\
&=-E(\alpha;\beta;\lambda/T)
E\left(\begin{array}{rr}\alpha+\gamma, & \beta+1 \\ \beta, & \beta+2 \end{array};-\frac{\lambda}{T}\right)+
\frac{1}{2\pi i}\int_{B}d\zeta\:
\frac{(\lambda/T)^{\zeta}\Gamma(\zeta)\Gamma(\alpha-\zeta)}{\Gamma(\beta-\zeta)} \\
&\times 
E\left(\begin{array}{rr}\alpha+\gamma, & \beta+l+1-\zeta \\ \beta, & \beta+l+2-\zeta \end{array};-\frac{\lambda}{T}\right),
\end{align*}

\noindent{}and therefore (\ref{eq:tl-1}) obeys the following Barnes integral representation

\begin{align}
\int_{0}^{T}t^{\beta+l}E(\alpha;\beta;\lambda/t)E(\alpha+\gamma;\beta;-\lambda/t)dt
=&\frac{T^{\beta+l+1}}{2\pi i} 
\int_{B}d\zeta\:\frac{(\lambda/T)^{\zeta}\Gamma(\zeta)\Gamma(\alpha-\zeta)}{\Gamma(\beta-\zeta)} \nonumber \\
&\times E\left(\begin{array}{rr}\alpha+\gamma, & \beta+l+1-\zeta \\ \beta, & \beta+l+2-\zeta \end{array};-\frac{\lambda}{T}\right).
\label{eq:tl-2}
\end{align}

\noindent{}By applying the recurrence relation [in \cite[\S 5.2.1, eq.~(9)]{Erd53} substitute $a_{1}=\alpha+\gamma$, 
$a_{2}=\beta+l+1-\zeta$, $\rho_{1}=\beta+1$, $\rho_{2}=\beta+l+2-\zeta$]

\begin{align*}
E\left(\begin{array}{rr}\alpha+\gamma, & \beta+l+1-\zeta \\ \beta, & \beta+l+2-\zeta \end{array};z\right)=&
\beta E\left(\begin{array}{rr}\alpha+\gamma, & \beta+l+1-\zeta \\ \beta+1, & \beta+l+2-\zeta \end{array};z\right) 
\nonumber \\
&-(1/z)
E\left(\begin{array}{rr}\alpha+\gamma+1, & \beta+l+2-\zeta \\ \beta+2, & \beta+l+3-\zeta \end{array};z\right)
\end{align*}

\noindent{}$l+1$ times, one derives

\begin{align}
E\left(\begin{array}{rr}\alpha+\gamma, & \beta+l+1-\zeta \\ \beta, & \beta+l+2-\zeta \end{array};z\right)=&
\sum_{n=0}^{l+1}\frac{(-1/z)^{n}\Gamma(\beta+l+1)}{\Gamma(\beta+n)}\binom{l+1}{n} \nonumber \\
&\times 
E\left(\begin{array}{rr}\alpha+\gamma+n, & \beta+l+n+1-\zeta \\ \beta+l+n+1, & \beta+l+n+2-\zeta \end{array};z\right).
\label{eq:Erecurl}
\end{align}

\noindent{}Substitute (\ref{eq:Erecurl}) in (\ref{eq:tl-2}), apply (\ref{eq:MacRob-2}) and (\ref{eq:hyper}), and get that

\begin{align}
\int_{0}^{T}t^{\beta+l}E(\alpha;\beta;\lambda/t)E(\alpha+\gamma;\beta;-\lambda/t)dt
=&
T^{\beta+l+1}
\sum_{n=0}^{l+1}\frac{(T/\lambda)^{n}\Gamma(\beta+l+1)}{\Gamma(\beta+n)}\binom{l+1}{n} \nonumber \\
&\times\Wh{\alpha-1}{\beta+l+n+1}{\gamma+n+1}{\beta}{\frac{\lambda}{T}}
\label{eq:Integrl2calc}
\end{align}

\noindent{}($\Re\beta+l>-1$). In virtue of (\ref{eq:MacRob-1}), combine (\ref{eq:Integrl2calc}) with
(\ref{eq:hyper3}), make a substitution $\lambda\to-1/\lambda$ and get (\ref{eq:integralMM}) as required.

To this end, we note that (\ref{eq:CoulombIntegral}) follows directly from (\ref{eq:integralMM}) and the definition of
the regular Coulomb wave function $F_{L}$.

\section{Properties of solution}

\subsection{Convergence}\label{sec:conv}

We prove that the series (\ref{eq:hyper}) converges uniformly for $\abs{z}\geq1$. Let

$$
\Wh{\alpha}{\beta}{\gamma}{\delta}{z}=
\frac{\Gamma(\alpha+1)}{\Gamma(\delta)}\sum_{\nu=0}^{\infty}u_{\nu}(z),\quad
u_{\nu}(z)\equiv 
\frac{z^{-\nu}\Gamma(\alpha+\gamma+\nu)}{\nu!\Gamma(\beta+1+\nu)}
\:_{2}F_{2}\left(\begin{array}{rr} \alpha+1, & \beta+\nu \\ \delta, & \beta+1+\nu \end{array};-\frac{1}{z}\right). 
$$

\noindent{}We show that for all $\varepsilon>0$, there exists $N=N(\varepsilon)$ such that, for $\nu>N$ and for 
$p\in\Z_{+}$,

\begin{equation}
\vert u_{\nu+1}(z)+u_{\nu+2}(z)+\ldots+u_{\nu+p}(z)\vert<\varepsilon\quad\text{for all }\quad\vert z\vert\geq1.
\label{eq:unu1}
\end{equation}

One can choose $N$ large enough so that for $\nu>N$, $\Gamma(\nu)\sim\sqrt{2\pi \nu}(\nu/e)^{\nu}$
(the Stirling's formula). But

$$
_{2}F_{2}\left(\begin{array}{rr} \alpha+1, & \beta+\nu \\ \delta, & \beta+1+\nu \end{array};-\frac{1}{z}\right)\sim 
M(\alpha+1;\delta;-1/z)
$$

\noindent{}for $\nu$ large, and thus

\begin{equation}
u_{\nu}(z)\sim \frac{e^{\beta+1-\alpha-\gamma}}{\sqrt{2\pi}}M(\alpha+1;\delta;-1/z)z^{-\nu}e^{\nu}\nu^{\alpha+\gamma-\beta-\nu-3/2}.
\label{eq:unu2}
\end{equation}

\noindent{}As seen, $\vert u_{\nu}(z)\vert\leq\vert u_{\nu}(1)\vert$ and
$u_{\nu}(z)\to 0$ as $\nu\to\infty$ for $\vert z\vert\geq1$. By applying the Minkowski's inequality to (\ref{eq:unu1}) and
replacing $\nu$ by $\nu+s$ in (\ref{eq:unu2}) for $s=1,2,\ldots,p$, one derives

\begin{align*}
& \Biggl\vert \sum_{s=1}^{p}u_{\nu+s}(z)\Biggr\vert\leq\sum_{s=1}^{p}\left\vert u_{\nu+s}(z)\right\vert\leq 
\sum_{s=1}^{p}\left\vert u_{\nu+s}(1)\right\vert \\
&\sim\frac{1}{\sqrt{2\pi}}\left\vert (e/\nu)^{\nu+\beta+1-\alpha-\gamma}\nu^{-1/2}M(\alpha+1;\delta;-1)\right\vert
\sum_{s=1}^{p}(e/\nu)^{s} \\
&\sim \frac{e}{\nu!}\left\vert (e/\nu)^{\beta+1-\alpha-\gamma}M(\alpha+1;\delta;-1)\right\vert<\varepsilon\quad
(\nu\;\text{large and}\;\vert z\vert\geq1)
\end{align*}

\noindent{}for all $p=1,2,\ldots$, where $\varepsilon>0$ can be chosen arbitrarily small. This proves that the series (\ref{eq:hyper})
converges uniformly for $\vert z\vert\geq1$.

Since $_{2}F_{2}(\cdot)$ is entire, that is, its radius of convergence is $\infty$, one deduces that the series is convergent
for $0<\vert z\vert<1$. We note that all other regions of convergence due to the parameters $\alpha,\beta,\gamma,\delta\in\C$ can be
obtained by applying the properties of $_{2}F_{2}(\cdot)$ to (\ref{eq:hyper}).

\subsection{Asymptotic expansion}\label{sec:asym}

The series expansion of (\ref{eq:hyper}) reads

$$
\frac{\Gamma(\alpha+1)\Gamma(\alpha+\gamma)}{\beta\Gamma(\beta)\Gamma(\delta)}+
\Biggl(\frac{\Gamma(\alpha+1)\Gamma(\alpha+\gamma+1)}{\Gamma(\beta+2)\Gamma(\delta)}
-\frac{\Gamma(\alpha+2)\Gamma(\alpha+\gamma)}{(\beta+1)\Gamma(\beta)\Gamma(\delta+1)}\Biggr)\cdot\frac{1}{z}
+O(z^{-2}).
$$

\noindent{}Hence,

\begin{equation}
\Wh{\alpha}{\beta}{\gamma}{\delta}{z}=\frac{\Gamma(\alpha+1)\Gamma(\alpha+\gamma)}{\Gamma(\beta+1)\Gamma(\delta)}+O(z^{-1})\quad
\text{as}\quad\vert z\vert\to\infty.
\label{eq:asymt}
\end{equation}

\noindent{}We note that the analysis of (\ref{eq:diffeq}) yields the same result. Indeed, functions $M(\cdot;\cdot;\pm1/z)\to1$
as $\vert z\vert\to\infty$, and so the right-hand side of (\ref{eq:diffeq}) becomes 
$[\Gamma(\alpha+1)\Gamma(\alpha+\gamma)/(\Gamma(\beta)\Gamma(\delta))]z^{-1}$ up to $O(z^{-2})$. By solving the modified 
differential equation one obtains (\ref{eq:asymt}).

By Paris \cite{Par05},

\begin{equation}
\:_{2}F_{2}\left(\begin{array}{rr} \alpha+1, & \beta+\nu \\ \delta, & \beta+1+\nu \end{array};-\frac{1}{z}\right)\sim 
(-1/z)^{\alpha-\delta}e^{-1/z}\frac{(\beta+\nu)\Gamma(\delta)}{\Gamma(\alpha+1)}\;\;\text{as}\;\;\vert z\vert\to0 
\label{eq:asym0}
\end{equation}

\noindent{}for $\vert\arg\:z\vert<\pi/2$ and $\Re z<0$. Substitute (\ref{eq:asym0}) in (\ref{eq:hyper}) and get

$$
\Wh{\alpha}{\beta}{\gamma}{\delta}{z}\sim (-1/z)^{\alpha-\delta}e^{-1/z}\sum_{\nu=0}^{\infty}
\frac{z^{-\nu}\Gamma(\alpha+\gamma+\nu)}{\nu!\Gamma(\beta+\nu)}
=\frac{\Gamma(\alpha+\gamma)}{\Gamma(\beta)}(-1/z)^{\alpha-\delta}e^{-1/z}M(\alpha+\gamma;\beta;1/z).
$$

\noindent{}But \cite[\S 13]{Abr72}
$M(\alpha+\gamma;\beta;1/z)\sim
(-z)^{\alpha+\gamma}\Gamma(\beta)/\Gamma(\beta-\alpha-\gamma)$ as $\vert z\vert\to0$ for $\Re z<0$,
and so 

\begin{equation}
\Wh{\alpha}{\beta}{\gamma}{\delta}{z}\sim
e^{-1/z}(-z)^{\gamma+\delta}\frac{\Gamma(\alpha+\gamma)}{\Gamma(\beta-\alpha-\gamma)}
\label{eq:asym01}
\end{equation}

\noindent{}as $\vert z\vert\to0$. Note that (\ref{eq:asym01}) as well as the asymptotic expansion of $f$ for $\Re z>0$
can be established from the direct inspection of (\ref{eq:diffeq}). In the latter case the right-hand side of (\ref{eq:diffeq}) obeys 
the form $e^{1/z}z^{\beta-\gamma}\Gamma(\alpha+1)/\Gamma(\delta-\alpha-1)$ ($\vert z\vert\to0;\Re z>0$), and thus

$$
\Wh{\alpha}{\beta}{\gamma}{\delta}{z}\sim
e^{1/z}z^{\beta+2-\gamma}\frac{\Gamma(\alpha+1)}{\Gamma(\delta-\alpha-1)}
$$

\noindent{}as $\vert z\vert\to0$ ($\Re z>0$).

\subsection{Some recurrence relations}\label{sec:recur}

By using the recurrence relations for $E$-functions [in \cite{Bho62}
substitute $p=q=2$, $\alpha_{1}=\alpha+1$, $\alpha_{2}=\beta+\nu$, $\rho_{1}=\delta$, $\rho_{2}=\beta+1+\nu$ and get (\ref{eq:recurr1});
in \cite[\S 5.2.1, eqs.~(7) and (9)]{Erd53} substitute $a_{1}=\alpha+1$, $a_{2}=\beta+\nu$, $\rho_{1}=\delta$, $\rho_{2}=\beta+1+\nu$ and 
get (\ref{eq:recurr2})--(\ref{eq:recurr3})] one obtains from (\ref{eq:hyper}) that

\begin{subequations}\label{eq:recurr}
\begin{align}
&E(\alpha+1;\delta;z)E(\alpha+\gamma;\beta+1;-z)=
\Wh{\alpha}{\beta}{\gamma}{\delta}{z}-\frac{1}{z}\Wh{\alpha+1}{\beta+1}{\gamma-1}{\delta+1}{z} \label{eq:recurr1} \\
&=\frac{\beta-\alpha-1}{\beta}
\Wh{\alpha}{\beta}{\gamma}{\delta}{z}+\frac{1}{\beta}\Wh{\alpha+1}{\beta}{\gamma-1}{\delta}{z}
+\frac{1}{\beta z^{2}}\Wh{\alpha+1}{\beta+2}{\gamma}{\delta+1}{z} \label{eq:recurr2} \\
&=\frac{\beta-\delta+1}{\beta}
\Wh{\alpha}{\beta}{\gamma}{\delta}{z}+\frac{1}{\beta}\Wh{\alpha}{\beta}{\gamma}{\delta-1}{z}
+\frac{1}{\beta z^{2}}\Wh{\alpha+1}{\beta+2}{\gamma}{\delta+1}{z}. \label{eq:recurr3}
\end{align}
\end{subequations}

\noindent{}Since the proofs are obvious due to the properties of $E$-functions we will omit them by noting that the recurrence relations
(\ref{eq:recurr2})--(\ref{eq:recurr3}) are reduced by applying (\ref{eq:recurr1}) so that the sum

$$
\sum_{\nu=0}^{\infty}\frac{z^{-\nu}\Gamma(\alpha+\gamma+\nu)}{\nu!\Gamma(\beta+\nu)}
E\left(\begin{array}{rr}\alpha+2, & \beta+1+\nu \\ \delta+1, & \beta+2+\nu \end{array};z\right)
$$

\noindent{}can be calculated by setting $\Gamma(\beta+\nu)=\Gamma(\beta+1+\nu)/(\beta+\nu)$, thus yielding

$$
\beta\Wh{\alpha+1}{\beta+1}{\gamma-1}{\delta+1}{z}+\frac{1}{z}\Wh{\alpha+1}{\beta+2}{\gamma}{\delta+1}{z}
$$

\noindent{}where the first summand is then represented by (\ref{eq:recurr1}).

\subsection*{Acknowledgments}

The author is very grateful to anonymous referee for pointing out Refs.~\cite{Hamza74,Shanker68} and
for many useful comments which have helped in improving the manuscript.


\section{Appendix : Alternative proof of Eq.~(\ref{eq:F2})}\label{sec:exp}

In \cite{Saad03} the authors have derived the integral (\ref{eq:F2}) 
[in \cite[Lemma~1]{Saad03} substitute $a=\alpha$, $b=b^{\prime}=\beta$, $a^{\prime}=
\alpha+\gamma$, $d=\beta+1$, $k=-k^{\prime}=i$]. Here, we shall demonstrate that 
(\ref{eq:F2}) is easy to obtain from (\ref{eq:tl-2}) by setting $\lambda=i$.

Recalling that $\exp(-ht)=\sum_{l=0}^{\infty}(-ht)^{l}/l!$, we can multiply (\ref{eq:tl-2}) by $(-h)^{l}/l!$ and then perform 
the summation over $l=0,1,\ldots$ This gives the integral

\begin{align*}
\int_{0}^{T}t^{\beta}e^{-ht}E(\alpha;\beta;i/t)E(\alpha+\gamma;\beta;-i/t)dt=&
\frac{T^{\beta+1}}{2\pi i}
\int_{B}d\zeta\:\frac{(i/T)^{\zeta}\Gamma(\zeta)\Gamma(\alpha-\zeta)}{\Gamma(\beta-\zeta)} \\
&\times\sum_{l=0}^{\infty}\frac{(-T h)^{l}}{l!}\:
E\left(\begin{array}{rr}\alpha+\gamma, & \beta+l+1-\zeta \\ \beta, & \beta+l+2-\zeta \end{array};-\frac{i}{T}\right).
\end{align*}

\noindent{}By applying (\ref{eq:MacRob-2}) to $E$-function,

\begin{align*}
&\int_{0}^{T}t^{\beta}e^{-ht}E(\alpha;\beta;i/t)E(\alpha+\gamma;\beta;-i/t)dt=
\frac{T^{\beta+1}}{2\pi i}
\int_{B}d\zeta\:\frac{(i/T)^{\zeta}\Gamma(\zeta)\Gamma(\alpha-\zeta)}{\Gamma(\beta-\zeta)} \nonumber \\
&\times\frac{1}{2\pi i}\int_{B}d\zeta^{\prime}\:\frac{(-i/T)^{\zeta^{\prime}}\Gamma(\zeta^{\prime})\Gamma(\alpha+\gamma-
\zeta^{\prime})}{\Gamma(\beta-\zeta^{\prime})}
\sum_{l=0}^{\infty}\frac{(-T h)^{l}}{l!}\cdot\frac{1}{\beta+l+1-\zeta-\zeta^{\prime}}. 
\end{align*}

\noindent{}The sum over $l$ equals $(T h)^{-\beta-1+\zeta+\zeta^{\prime}}(\Gamma(\beta+1-\zeta-\zeta^{\prime})-
\Gamma(\beta+1-\zeta-\zeta^{\prime},T h))$, where the incomplete gamma function 
$\Gamma(a,t)=\int_{t}^{\infty}u^{a-1}e^{-u}du$.

By the residue theorem,

\begin{align}
&\int_{0}^{T}t^{\beta}e^{-ht}E(\alpha;\beta;i/t)E(\alpha+\gamma;\beta;-i/t)dt=
\frac{\beta\Gamma(\alpha)\Gamma(\alpha+\gamma)}{h^{\beta+1}\Gamma(\beta)} \nonumber \\
&\times F_{2}(\beta+1;\alpha,\alpha+\gamma;\beta,\beta;i/h,-i/h)-\frac{h^{-\beta-1}}{2\pi i}\int_{B}d\zeta\:
\frac{(ih)^{\zeta}\Gamma(\zeta)\Gamma(\alpha-\zeta)}{\Gamma(\beta-\zeta)} \nonumber \\
&\times\frac{1}{2\pi i}
\int_{B}d\zeta^{\prime}\:\frac{(-ih)^{\zeta^{\prime}}\Gamma(\zeta^{\prime})\Gamma(\alpha+\gamma-
\zeta^{\prime})\Gamma(\beta+1-\zeta-\zeta^{\prime},T h)}{\Gamma(\beta-\zeta^{\prime})}
\quad(\abs{h}\geq2), 
\label{eq:test1}
\end{align}

\noindent{}where the contour integral involving the incomplete gamma function is calculated by using the integral representation of
$\Gamma(\cdot,T h)$ along with the residue theorem, namely,

\begin{align}
&\frac{1}{2\pi i}
\int_{B}d\zeta^{\prime}\:\frac{(-ih)^{\zeta^{\prime}}\Gamma(\zeta^{\prime})\Gamma(\alpha+\gamma-
\zeta^{\prime})\Gamma(\beta+1-\zeta-\zeta^{\prime}, T h)}{\Gamma(\beta-\zeta^{\prime})} \nonumber \\
&=\sum_{\nu=0}^{\infty}\frac{(-i/h)^{\nu}\Gamma(\alpha+\gamma+\nu)\Gamma(\beta+1-\zeta+\nu, T h)}{\nu!\Gamma(\beta+\nu)}=
\int_{ T h}^{\infty}dt\:t^{\beta-\zeta}e^{-t}
\sum_{\nu=0}^{\infty}\frac{(-it/h)^{\nu}\Gamma(\alpha+\gamma+\nu)}{\nu!\Gamma(\beta+\nu)} \nonumber \\
&=
\int_{ T h}^{\infty}t^{\beta-\zeta}e^{-t}E(\alpha+\gamma;\beta;-ih/t)dt
=h^{\beta+1-\zeta}\int_{ T}^{\infty}t^{\beta-\zeta}e^{-ht}E(\alpha+\gamma;\beta;-i/t)dt, \label{eq:test2}
\end{align}

\noindent{}where in the last step the substitution $t\to ht$ has been initiated.

Substitute (\ref{eq:test2}) in (\ref{eq:test1}) and get by (\ref{eq:MacRob-2}),

\begin{align*}
\int_{0}^{ T}t^{\beta}e^{-ht}E(\alpha;\beta;i/t)E(\alpha+\gamma;\beta;-i/t)dt=&
\frac{\beta\Gamma(\alpha)\Gamma(\alpha+\gamma)}{h^{\beta+1}\Gamma(\beta)}
\:F_{2}(\beta+1;\alpha,\alpha+\gamma;\beta,\beta;i/h,-i/h) \\
&-\int_{ T}^{\infty}t^{\beta}e^{-ht}E(\alpha;\beta;i/t)E(\alpha+\gamma;\beta;-i/t)dt.
\end{align*}

\noindent{}But $\int_{0}^{ T}+\int_{ T}^{\infty}=\int_{0}^{\infty}$ for all $ T\in\R$, and hence (\ref{eq:F2}) holds.

\section*{References}

\bibliographystyle{unsrt}

\newcommand{\nosort}[1]{}

\end{document}